\documentclass[12pt]{article}
\usepackage{amsthm,amsmath,amssymb,amscd,verbatim}
\usepackage{graphics,pifont}
\usepackage{amssymb}
\usepackage{graphicx}
\date{}

\newcommand{\va}{\varepsilon}
\newcommand{\sq}{$\square$}

\begin{document}
\title{What make them all so turbulent}    
\author{Bau-Sen Du \\  
Institute of Mathematics \cr
Academia Sinica \cr
Taipei 10617, Taiwan \cr
dubs@math.sinica.edu.tw \cr}
\maketitle

\begin{abstract}
\noindent
We give a unified proof of the existence of turbulence for some classes of continuous interval maps which include, among other things, maps with periodic points of odd periods $> 1$, some maps with dense chain recurrent points and densely chaotic maps. 


\end{abstract}

Let $I$ be a compact interval in the real line and let $f : I \to I$ be a continuous map.  It is well-known that {\bf\cite{ba, bc, bl, du, ru}} if (a) there exist a point $c$ and an odd integer $n > 1$ such that $f^n(c) \le c < f(c)$ or $f(c) < c \le f^n(c)$, or (b) $f$ has dense periodic points and $f^2(a) \ne a$ for some point $a$, or (c) there is a point whose $\omega$-limit set with respect to $f$ contains a fixed point $z$ of $f$ and a point $\ne z$, or (d) $f$ is densely chaotic, i.e., the set $LY(f) = \{ (x, y) \in I \times I : \limsup_{n \to \infty} |f^n(x) - f^n(y)| > 0$ and $\liminf_{n \to \infty} |f^n(x) - f^n(y)| = 0 \}$ is dense in $I \times I$, then $f^2$ is turbulent (and $f$ has periodic points of all even periods).  Since turbulent maps are known {\bf\cite{bc}} to be topologically semi-conjugate, on some compact invariant subsets, to the shift map on two symbols which is a typical model for chaotic dynamical systems, these maps $f^2$ (and so $f$) are chaotic.  When we examine closely the above 4 conditions, we find that none is implied by {\it all} other three (see Figures 1 \& 2).  So, what do they have in common which make them all so turbulent?  In this note, we answer this question by a simple result (Theorem 1) which extends Proposition 3 on page 122 of {\bf\cite{bc}}.

\begin{figure}[ht] 
\begin{center}
\includegraphics[width=6cm,height=5cm]{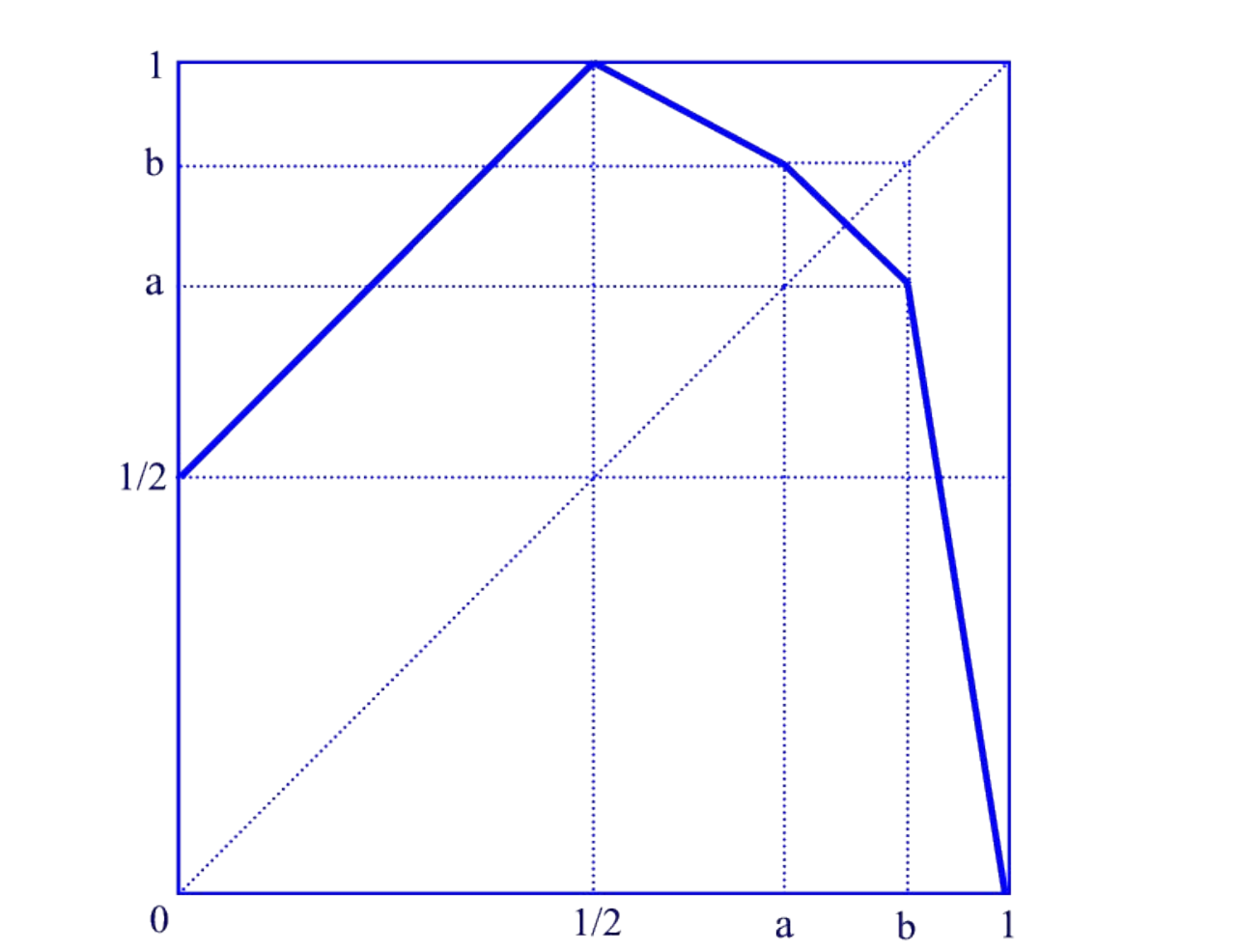} 
\caption{A map satisfying (a), but none of (b), (c) and (d).}
\end{center}
\end{figure}
\begin{figure}[ht] 
\begin{center}
\includegraphics[width=6cm,height=5cm]{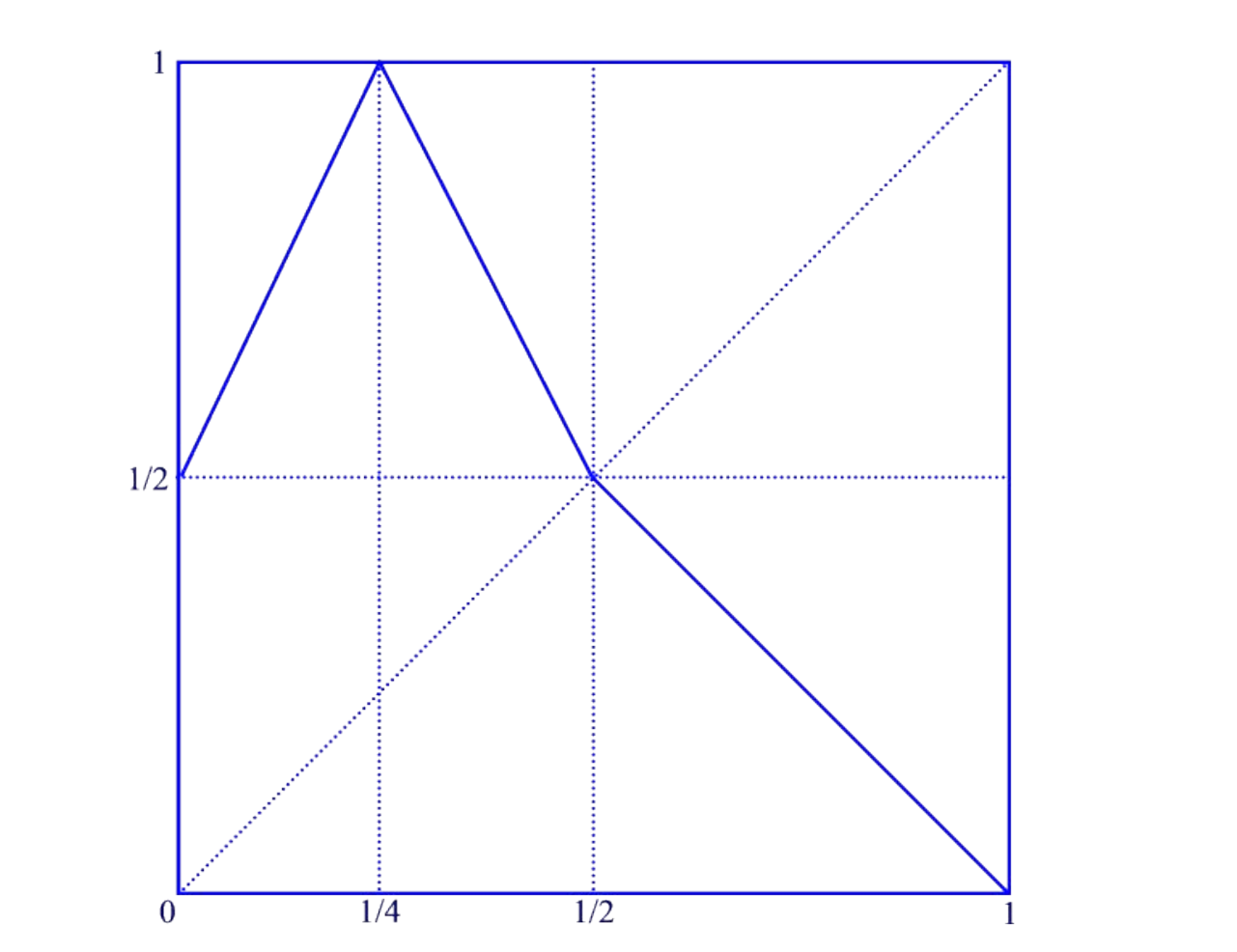} 
\caption{A map satisfying (b), (c) and (d), but not (a).}
\end{center}
\end{figure}

Let $J$ be a compact interval in $I$.  If there exist two compact subintervals $J_0$ and $J_1$ of $J$ with at most one point in common such that $f(J_0) \cap f(J_1) \supset J_0 \cup J_1$, then we say that $f$ is turbulent on $J$ (and on $I$) {\bf\cite{bc}}.  If there exist two compact subintervals $K$ and $L$ of $I$ with at most one point in common such that $f$ is turbulent on $K$ and on $L$, then we say that $f$ is doubly turbulent on $I$.  

\noindent
{\bf Theorem 1.} 
{\it Let $f$ be a continuous map from $I$ into itself and let $x_0$ be a point in $I$.  Then exactly one of the following holds:
\begin{itemize}
\item[(A)]
If there exist a point $c$ in the orbit $O_f(x_0) = \{ x_0, f(x_0), f^2(x_0), \cdots \}$ of $x_0$ and an integer $n \ge 2$ such that $f^n(c) \le c < f(c)$ or $f(c) < c \le f^n(c)$, then at least one of the following holds:
\begin{itemize}
\item[(1)]
There exist a fixed point $z$ of $f$ and a compact subinterval $K$ of $I$ such that (i) $c \in K$, (ii) $f^2(K) \subsetneq K$, (iii) $K$ contains no fixed points of $f$, and (iv) $K$ and $f(K)$ lie on opposite sides of $z$, in particular, the iterates of $c$ with respect to $f$ are "`jumping"' alternately around the fixed point $z$;

\item[(2)]
$f$ has periodic points of all even periods and $f^2$ is doubly turbulent.
\end{itemize}

\item[(B)]
If $x_i = f^i(x_0)$ for all $i \ge 0$, then either for some $m > 0$ the sequence $< x_n >_{n \ge m}$ converges monotonically to a fixed point of $f$ or there exist a fixed point $\hat z$ of $f$ and a strictly increasing sequence $0 \le n_0 < n_1 < n_2 < \cdots$ of integers such that if $x_0 < x_1$ (if $x_0 > x_1$ then all inequalities below are reversed) then 
\begin{small}
$$x_0 < x_1 < \cdots < x_{n_0-1} \quad < x_{n_1} < x_{n_1+1} < \cdots < x_{n_2-1} \quad < x_{n_3} < x_{n_3+1} < \cdots < x_{n_4-1} \quad < \cdots < \hat z$$ $$ < \cdots < x_{n_5-1} < \cdots < x_{n_4+1} < x_{n_4} \quad < x_{n_3-1} < \cdots < x_{n_2+1} < x_{n_2} \quad < x_{n_1-1} < \cdots < x_{n_0+1} < x_{n_0}$$ and if $p = \lim_{i \to \infty} x_{n_{2i+1}}$ and $q = \lim_{i \to \infty} x_{n_{2i}}$ then $p \le \hat z \le q$ and, $f(p) = q$ and $f(q) = p$.  
\end{small}  
In particular, $x_0$ is asymptotically periodic of period 1 or 2, i.e., there is a periodic point $y$ of $f$ with $f^2(y) = y$ such that $\lim_{n \to \infty} |f^n(x_0) - f^n(y)| = 0$.  
\end{itemize}}

\noindent
{\it Proof.}
If the hypothesis of {\it (A)} fails, then it is clear that {\it (B)} holds.  Now, assume that $f^n(c) \le c < f(c)$ for some point $c$ in $O_f(x_0)$.  If $f(c) < c \le f^n(c)$, the proof is similar.  Let $X = \{ f^i(c) : 0 \le i \le n-1 \}$.  Let $a = \max \{ x \in X : f(x) > x \}$ and let $b$ be any point in $X \cap [a, f(a)]$ such that $f(b) \le a$.  Then $c \le a$.  Let $z$ be a fixed point of $f$ in $[a, b]$ and let $v$ be a point in $[a, z]$ such that $f(v) = b$.  So, $f^2(v) = f(b)$ and $\max \{ c, f^2(v) \} \le a \le v < z < b = f(v)$.  Let $z_0 = \min \{ v \le x \le z : f^2(x) = x \}$.  Then $f(x) > z$ and $f^2(x) < x$ for all $v \le x < z_0$.  We have three cases to consider:

Case 1. If $f^2(x) < z_0$ for all $\min I \le x \le v$, then $f^2(x) < z_0 \le z < f(x)$ for all $\min I \le x \le z_0$.  Let $t$ be a point in $(v, z_0)$ such that $t >  f^2(x)$ for all $\min I \le x \le t$.  Let $K = [\min I, t]$.  Then $c \in [\min I, v] \subset K$, $f^2(K) \subset [\min I, t) \subsetneq K$, $K$ contains no fixed points of $f$, and $K$ and $f(K)$ lie on opposite sides of $z$.  

Case 2. If the point $d = \max \{ \min I \le x \le v : f^2(x) = z_0 \}$ exists and $\min \{ f^2(x) : d \le x \le z_0 \}$ $= s > d$, then $f(x) > z \ge z_0 > f^2(x) \ge s$ for all $d < x < z_0$.  Let $\tilde t$ be a point in $(v, z_0)$ such that $\tilde t > f^2(x)$ for all $s \le x \le v$.  Let $K = [s, \tilde t]$.  Then $K$ contains no fixed points of $f$, $K$ and $f(K)$ lie on opposite sides of $z$ and $f^2(K) \subset [s, \tilde t) \subsetneq K$.  Furthermore, for some $2 \le k \le n$, $f^{k-1}(c) = b$ and so, $f^k(c) = f(b) = f^2(v) \in f^2(K) \subset K$.  Consequently, $f^k(c) \in K$.  Since $f(K \cup f(K)) \subset K \cup f(K)$ and $n \ge k$, we have $f^n(c) = f^{n-k}(f^k(c)) \in K \cup f(K)$. Since $f^n(c) \, (\le c \le v) < z$, this forces $f^n(c) \in K$.  Since $\tilde t \in [s, \tilde t] = K$ and $f^n(c) \le c \le v < \tilde t$, this in turn implies that $c \in K$.    

Case 3. If both the point $d = \max \{ \min I \le x \le v : f^2(x) = z_0 \}$ and the point $u_1 = \min \{ d \le x \le z_0 : f^2(x) = d \}$ exist, then $f(x) > z \ge z_0 > f^2(x)$ on $(d, z_0)$ and $f^2([d, u_1]) \cap f^2([u_1, z_0]) \supset [d, z_0] = [d, u_1] \cup [u_1, z_0]$.  In particular, $f^2$ is turbulent on $[d, z_0] \subset [\min I, z]$.  Furthermore, since $u_1 = \min \{ d \le x \le z_0 : f^2(x) = d \}$, we have $d < f^2(x) < z_0$ on $(d, u_1)$.  Let $p_1$ be any point in $(d, u_1)$ such that $f^2(p_1) = p_1$.  Let $u_2 = \min \{ d \le x \le p_1 : f^2(x) = u_1 \}$.  Then $d < (f^2)^2(x) < z_0$ on $(d, u_2)$.  Let $p_2$ be any point in $(d, u_2)$ such that $(f^2)^2(p_2) = p_2$.  Inductively, we obtain points $d < \cdots < p_n < u_n < \cdots < p_2 < u_2 < p_1 < u_1 < z_0$ such that $u_n = \min \{ d \le x \le p_{n-1} : (f^2)^{n-1}(x) = u_1 \}$, $d < (f^2)^n(x) < z_0$ on $(d, u_n)$ and $(f^2)^n(p_n) = p_n$.  Since $f(x) > z \ge z_0$ on $(d, z_0)$, we have $f^i(p_n) < z_0 < f^j(p_n)$ for all even $i$ and all odd $j$ in $[0, 2n]$.  So, each $p_n$ is a period-$(2n)$ point of $f$.  This confirms that $f$ has periodic points of all {\it even} periods.  Finally, since $d$ is the largest point in $[\min I, z_0)$ such that $f^2(d) = z_0$, $f$ must map the endpoints of $[d, z_0]$ {\it into} the endpoints of $f([d, z_0])$ and no points $x$ in $(d, z_0)$ can satisfy $f(x) = f(d)$ or $f(x) = f(z_0)$.  Consequently, if $f(d) > f(z_0)$ (if $f(d) < f(z_0)$, the proof is similar), then $f([d, z_0]) = [f(z_0), f(d)]$ and, for some $\hat s \le d$, $f((f(z_0), f(d)) = f^2((d, z_0)) = [\hat s, z_0) \supset [d, z_0)$.  Let $e$ be a point in $(f(z_0), f(d))$ such that $f(e) = d$.  Then $f^2([f(z_0), e]) \cap f^2([e, f(d)]) \supset [f(z_0), f(d)] = [f(z_0), e] \cup [e, f(d)]$.  Furthermore, if $f(d) = f(z_0)$ and $f([d, z_0]) = [r, f(d)]$ for some point $r > z$ (if $f(d) = f(z_0)$ and $f([d, z_0]) = [f(d), r]$, the proof is similar), then since $f([r, f(d)]) = f^2([d, z_0]) \supset [d, z_0]$, there exists a point $u$ in $[r, f(d))$ such that $f(u) = d$.  Since $f^2([u, f(d)]) \supset f([d, z_0]) = [r, f(d)]$, there exists a point $w$ in $(u, f(d))$ such that $f^2(w) = r$.  Therefore, $f^2([u, w]) \cap f^2([w, f(d)]) \supset [r, f(d)] \supset [u, f(d)] = [u, w] \cup [w, f(d)]$.  In either case, $f^2$ is turbulent on $[z, \max I]$.  This, combined with the above, shows that $f^2$ is doubly turbulent on $I$.
\hfill\sq

In Part{\it (A)(1)} of the above result, the compact interval $K$ is not an ordinary one.  It is one with the following 4 properties that (i) $c \in K$; (ii) $f^2(K) \subsetneq K$; (iii) $K$ contains no fixed points of $f$; and (iv) $f(K) \cap K = \emptyset$.  By choosing the appropriate point $c$, it is the violation of one of these properties that establishes the following result in which {\it (2)} and {\it (4)} are generalizations of (b) and (d) above respectively.  

\noindent
{\bf Corollary 2.} 
{\it Each of the following statements implies that $f$ has periodic points of all even periods and $f^2$ is doubly turbulent:  
\begin{itemize}
\item[(1)]
There exist a point $c$ and an odd integer $n > 1$ such that $f^n(c) \le c < f(c)$ or $f(c) < c \le f^n(c)$, in particular, $f$ has a periodic point of odd period $> 1$;

\item[(2)]
The chain recurrent points of $f$ are dense in $I$ and $f^2(a) \ne a$ for some point $a$ in $I$ (recall that a chain recurrent point is a point $x$ which satisfies that for every $\va > 0$ there exist a finite sequence of points $x_i, 0 \le i \le n$ such that $x_0 = x = x_n$ and $|f(x_i) - x_{i+1}| < \va$ for all $0 \le i \le n-1$.  
Note that if $x_0$ is a chain recurrent point of $f$ with $f(x_0) < x_0$ \, ($f(x_0) > x_0$ respectively), then by discussing the three cases similar to those three in the above proof of Theorem 1 with $x_0$ replacing $z_0$, we can obtain (see Lemma 32 on page 150 of {\bf\cite{bc}}) a point $c$ such that $f(c) < c < f^2(c) = x_0$ \, ($x_0 = f^2(c) > c > f(c)$ respectively));

\item[(3)]
The $\omega$-limit set $\omega_f(b)$ of some point $b$ in $I$ contains a fixed point $z$ of $f$ and a point $\ne z$;

\item[(4)]
There is a point in $I$ which is not asymptotically periodic of period 1 or 2 and the set $\{ (x, y) : \liminf_{n \to \infty} |f^n(x) - f^n(y)| = 0 \}$ is dense in $I \times I$, in particular, $f$ is densely chaotic; 

\item[(5)]
There is a point $c$ in $I$ such that $$\limsup_{n \to \infty} |f^n(c) - f^{n+1}(c)| > 0 \quad \text{and} \quad \liminf_{n \to \infty} |f^n(c) - f^{n+1}(c)| = 0.$$
\end{itemize}}

The following result can be proved similarly.

\noindent
{\bf Theorem 3.}
{\it If there exist a fixed point $z$ of $f$ and a point $c$ of $I$ such that $f(c) < c < z$ or $z < c < f(c)$, then at least one of the following holds:
\begin{itemize}
\item[(1)]
$f$ has a proper compact interval $J$ in $I$ such that $c \in J$, $f(J) \subsetneq J$ and $z \notin J$;

\item[(2)]
$f$ is turbulent and has periodic points of all periods.
\end{itemize}

\noindent
Consequently, if 
(1) there exist a fixed point $z$ of $f$, a point $c$ of $I$ and an integer $n \ge 2$ such that $f(c) < c < z \le f^n(c)$ or, $f^n(c) \le z < c < f(c)$, or $(f(c) - z)/(c - z) > 1$ and $z \in \omega_f(c)$; or 
(2) the chain recurrent points of $f$ are dense in $I$ and $f$ has at least two fixed points and $f(a) \ne a$ for some point $a$ in $I$, then $f$ is turbulent and has periodic points of all periods.}

\end{document}